\documentclass[aos,MSNbibl,nameyear,dvips]{arximspdf}
\usepackage{graphicx}

\doi{10.1214/12-AOS991} 
\volume{40}
\issue{2}
\pubyear{2012}
\firstpage{1061}
\lastpage{1073}

\makeatletter
\newtheorem{them}{Theorem}
\newproclaim{definition}{Definition}
\newtheorem{proposition}{Proposition}
\newtheorem{corollary}{Corollary}

\newcommand\cip{\protect\mathpalette{\protect\independenT}{\perp}}
\def\independenT#1#2{\mathrel{\rlap{$#1#2$}\mkern4.1mu{#1#2}}}
\newcommand{\trace}{\operatorname{tr}}
\def\bsuffix #1{#1}
\makeatother

\begin{document}
\begin{frontmatter}

\title{Estimation of means in graphical Gaussian models with symmetries}
\runtitle{Estimation of means}

\begin{aug}
\author{\fnms{Helene} \snm{Gehrmann}\ead[label=e1]{gehrmann@stats.ox.ac.uk}}
\and
\author{\fnms{Steffen L.} \snm{Lauritzen}\corref{}\ead[label=e2]{steffen@stats.ox.ac.uk}}
\runauthor{H. Gehrmann and S. L. Lauritzen}
\affiliation{University of Oxford}
\address{Department of Statistics\\
University of Oxford\\
1 South Parks Road\\
Oxford OX1 3TG\\
United Kingdom\\
\printead{e1}\\
\hphantom{E-mail: }\printead*{e2}} 
\end{aug}

\received{\smonth{2} \syear{2011}}
\revised{\smonth{10} \syear{2011}}

%
\begin{abstract}
We study the problem of estimability of means in undirected graphical
Gaussian models with symmetry restrictions represented by a colored
graph. Following on from previous studies, we partition the variables
into sets of vertices whose corresponding means are restricted to being
identical. We find a necessary and sufficient condition on the
partition to ensure equality between the maximum likelihood and
least-squares estimators of the mean.
\end{abstract}

%

\end{frontmatter}


\section{Introduction}

The elegant principles of symmetry and invariance appear in many areas
of statistical research [e.g., \citet{dawid}, \citet{diaconis},
\citet{eaton2}, \citet{viana}].
Symmetry restrictions in the multivariate Gaussian distribution have a
long history [\citet{andersson}, \citet{anderssonbrons},
\citet{jensen}, \citet{olkin72}, \citet{olkin}, \citet{votaw},
 \citet{wilks}] and have recently been combined with
conditional independence relations [\citet{andersen}, \citet{lauritzensym}, \citet{hylleberg},
 \citet{madsen}].

This article is concerned with graphical Gaussian models with symmetry
constraints introduced by \citet{lauritzensym}. The types of
restrictions are: equality between specified elements of the
concentration matrix (RCON), equality between specified partial
correlations (RCOR) and restrictions generated by permutation symmetry
(RCOP), a special instance of the former two. The models can be
represented by vertex and edge colored graphs, where parameters
associated with equally colored vertices or edges are restricted to
being identical.

We consider maximum likelihood estimation of a nonzero mean $\mu$
subject to specific equality constraints, assuming the covariance
matrix $\Sigma$ satisfies the restrictions of one of the above models.
This could be relevant, for example, if treatment effects are to be
estimated in an experiment where the basic error structure in the units
exhibits conditional independencies in a~symmetric pattern.

For the Gaussian distribution, maximum likelihood estimation of $\mu$
under an unknown covariance structure is generally nontrivial, as the
maximizer of the likelihood function in $\mu$ for fixed $\Sigma$ may
depend on $\Sigma$. The least-squares estimator $\mu^{*}$, however,
is defined by minimizing the sum of squares so that in case of equality
of $\hat{\mu}$ and $\mu^{*}$, the former is independent of $\Sigma
$. \citet{kruskal} showed that for for fixed $\Sigma$, and $\mu$ in
an affine space $\Omega$, $\hat{\mu}$ and $\mu^{*}$ agree if and
only if $\Omega$ is stable under $\Sigma$, or equivalently under $K =
\Sigma^{-1}$; see also \citet{haberman} and \citet{eaton}. Here we
derive a necessary and sufficient condition on the graph coloring
representing a model and the symmetry constraints on the mean vector
$\mu$ which ensures this stability and hence equality of estimators.

We let $G=(V,E)$ denote the dependence graph of the model and let its
colored version be denoted by $\mathcal{G}=(\mathcal{V},\mathcal
{E})$, where $\mathcal{V}$ denotes a partition of $V$ into vertex
color classes and $\mathcal{E}$ a partition of $E$ into edge color
classes. The symbol $\mathcal{M}$ is to denote a partition of $V$ such
that whenever two vertices are in the same set of $\mathcal{M}$, the
corresponding means are restricted to being equal. The necessary and
sufficient condition for equality of $\hat{\mu}$ and $\mu^{*}$ in
the symmetry model represented by $(\mathcal{V},\mathcal{E})$ is
twofold: (i) the partition $\mathcal{M}$ must be finer than $\mathcal
{V}$; and (ii) the partition must be \emph{equitable} with respect to
every edge color class in $\mathcal{E}$ as defined by \citet{sachs}.

\begin{figure}

\includegraphics{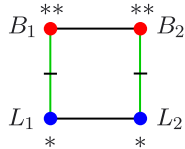}

\caption{Graphical Gaussian symmetry model supported by Frets's heads.}
\label{fig:frets}
\end{figure}

For example, the graph in Figure~\ref{fig:frets} represents a
graphical Gaussian symmetry model for data concerning the head
dimensions of first and second sons known as \emph{Frets's heads}
[\citet{frets}, \citet{mardia}]; here $L_1, B_1$ denotes the length and breadth
of the head of the first son, and similarly for $L_2,B_2$.

This model has previously been found to be well supported by the data
[\citet{gehrmann}, \citet{lauritzensym}, \citet{whittaker}] when no constraints were
considered on the means. We may, for example, be interested in the
hypothesis that the two lengths have the same mean, and the two
breadths have the same mean, indicating that head dimensions do not
generally change with the parity of the son. We shall demonstrate that
this hypothesis is simple in the sense that the maximum likelihood
estimator, or MLE for short, of the mean under this hypothesis can be
found by a simple average. Also, we shall demonstrate that this is not
the case if we consider lengths and breadths separately.
\section{Preliminaries and notation}
\subsection{Graphical Gaussian models}\label{ugm}
Let $G=(V,E)$ be an undirected graph with vertex set $V$ and edge set
$E$ and let $Y=(Y_{\alpha})_{\alpha\in V}$ be a multivariate Gaussian
random vector. Then the \emph{graphical Gaussian model} represented by
$G$ is the set of Gaussian distributions for which $Y_{\alpha}$ is
conditionally independent of $Y_{\beta}$ given the remaining
variables, denoted $Y_{\alpha} \cip Y_{\beta} \mid Y_{V \setminus\{
\alpha,\beta\}}$, whenever there is no edge between $\alpha$ and
$\beta$ in $G$.

For a multivariate Gaussian $\mathcal{N}_{|V|}(\mu, \Sigma)$
distribution with concentration matrix $\Sigma^{-1} = K = (k_{\alpha
\beta})_{\alpha, \beta\in V}$, it holds that $Y_{\alpha} \cip
Y_{\beta} \mid Y_{V \setminus\{\alpha,\beta\}}$ if and only if
$k_{\alpha\beta} = 0$. Thus letting $\mathcal{S}^{+}(G)$ denote the
set of symmetric positive definite matrices indexed by $V$ whose
$\alpha\beta$-entry is zero whenever $\alpha\beta\notin E$, the
graphical Gaussian model represented by $G$ assumes
\[
Y \sim\mathcal{N}_{|V|}(\mu, \Sigma),\qquad \mu\in\Omega=\mathbb
{R}^V, \Sigma^{-1} = K \in\mathcal{S}^{+}(G).
\]
For further details, see, for example, \citet{lauritzenbook}, Chapter~5.

\subsection{Graph coloring}
For general graph terminology we refer to \citet{bollobasbook}.
Following \citet{lauritzensym} we define the following notation for
graph colorings. Let $G=(V,E)$ be an undirected graph. Then a \emph
{vertex coloring} of $G$ is a partition $\mathcal{V}=\{V_{1}, \ldots,
V_{k}\}$ of $V$, where we refer to $V_{1}, \ldots, V_{k}$ as the \emph
{vertex color classes}. Similarly, an \emph{edge coloring} of $G$ is a~partition $\mathcal{E}=\{E_{1}, \ldots, E_{l}\}$ of $E$ into $l$
\emph{edge color classes} $E_{1}, \ldots, E_{l}$. We call a color
class containing one element a \emph{singleton} and a partition
containing only singletons as elements a \emph{singleton partition}.

Then $\mathcal{G}=(\mathcal{V},\mathcal{E})$ denotes the \emph
{colored graph} with vertex coloring $\mathcal{V}$ and edge coloring
$\mathcal{E}$; we also say that $(\mathcal{V},\mathcal{E})$ is a
\emph{graph coloring}. We indicate the color class of a vertex by the
number of asterisks we place next to it. Similarly we indicate the
color class of an edge by dashes. color classes which are singletons
are displayed in black and without asterisks or dashes.

\subsection{Further notation}\label{section_notation}
As we shall be considering constraints on the mean vector defined by
partitions of the mean vector into groups of equal entries, we
introduce the following notation. For $\mathcal{M}$ a partition of $V$
and $\alpha, \beta\in V$, we write $\alpha\equiv\beta (\mathcal
{M})$ to denote that $\alpha$ and $\beta$ lie in the same set in~$\mathcal{M}$ and let $\Omega(\mathcal{M})$ be the linear space of
vectors which are constant on each set of the partition $\mathcal{M}$:
%
\begin{equation}\label{omega_m}
\Omega(\mathcal{M}) = \{(x_{\alpha})_{\alpha\in V} \in\mathbb
{R}^{V}\dvtx x_{\alpha} = x_{\beta} \mbox{ whenever } \alpha\equiv
\beta (\mathcal{M})\}.
\end{equation}

For two partitions $\mathcal{M}_{1}$ and $\mathcal{M}_{2}$ of the
same set, we shall say that $\mathcal{M}_{1}$ is \emph{finer than}
$\mathcal{M}_{2}$, denoted by $\mathcal{M}_{1} \leq\mathcal
{M}_{2}$, if every set in $\mathcal{M}_{2}$ can be expressed as a
union of sets in $\mathcal{M}_{1}$. We equivalently say that $\mathcal
{M}_{2}$ is \emph{coarser than} $\mathcal{M}_{1}$.

If $A$ is a set of edges in a graph $G$, for $\alpha\in V$ we shall
write $\operatorname{ne}_{A}(\alpha)$ to denote the set of vertices
which are connected to $\alpha$ by an edge inside $A$.

We further adopt the following notation from \citet{lauritzensym}.
For a colored graph $\mathcal{G}=(\mathcal{V},\mathcal{E})$ and $u
\in\mathcal{V}$ we let $T^{u}$ denote the $|V| \times|V|$ diagonal
matrix with $T^{u}_{\alpha\alpha} = 1$ if $\alpha\in u$ and 0
otherwise. Similarly, for each edge color class $u \in\mathcal{E}$ we
let $T^{u}$ be the $|V| \times|V|$ symmetric matrix with
$T^{u}_{\alpha\beta} = 1$ if $\{\alpha,\beta\} \in u$ and 0 otherwise.

\section{Maximum likelihood and least-squares estimation}\label{mle_section}
Letting $K = \Sigma^{-1}$ as above, the density of $Y \sim\mathcal
{N}_{|V|}(\mu,\Sigma)$ is given by
\[
f_{\mu,K}(y)=\frac{\det K^{1/2}}{(2\pi)^{{|V|/2}}}\exp\{-(y-\mu
)^{T}K(y-\mu)/2 \}
\]
so that the likelihood function based on a sample $\mathcal{Y} =
(Y^{i})_{1 \leq i \leq n}$ where $Y^i$ are independent, and $Y^i \sim
\mathcal{N}_{|V|}(\mu,\Sigma)$ becomes
%
\begin{equation}\label{eq:likfunction}
L(\mu,K;y)\propto\det K^{n/2}\exp\biggl\{-\sum_{1\leq i \leq n}(y^i-\mu
)^{T}K(y^i-\mu)/2 \biggr\}.
\end{equation}

If $\mu$ is unrestricted, so that $\mu\in\Omega= \mathbb{R}^V$,
the likelihood function $L$ in (\ref{eq:likfunction}) is maximized
over $\mu$ for fixed $K$ by the least squares estimator $\mu^*=\bar
{y}$, and inference about $K$ can be based on the profile likelihood function
%
\begin{equation}\label{eq:proflik}
L(\hat\mu,K;y)\propto\det K^{n/2}\exp\{-\trace(KW)/2\},
\end{equation}
where $W=\sum_{i=1}^{n}(y^{i}-\mu^*)(y^{i}-\mu^*)^{T}$ is the matrix
of sums of squares and products of the residuals.
However, inference about $\mu$ when $K$ is unknown and needs to be
estimated is generally not possible, a classic example being known as
the \emph{Behrens--Fisher problem} [see \citet{scheffe} and \citet
{drtonbehrens}], where $\Sigma$ is bivariate and diagonal whereas the
mean vector satisfies the restriction $\mu_1=\mu_2$.

\citet{kruskal} found the following necessary and sufficient condition
for the two estimators to agree for fixed $\Sigma$:
\begin{them}[(Kruskal)]\label{thm:kruskal}
Let $Y$ be a random vector in an inner product space with unknown mean
$\mu$ in a linear space $\Omega$ and known covariance matrix~$\Sigma
$. Then the estimators $\mu^{*}$ and $\hat{\mu}$ coincide if and
only if $\Omega$ is invariant under $K=\Sigma^{-1}$, that is, if and
only if
%
\begin{equation}\label{eq_kruskal}
K\Omega\subseteq\Omega.
\end{equation}
\end{them}

Consequently, if (\ref{eq_kruskal}) is satisfied by all $K$ in a
model, the likelihood function can be maximized for fixed $K$ by $\mu
^*$, and inference on $K$ can be based on the profile likelihood (\ref
{eq:proflik}).

\section{Model types}
\subsection{Descriptions}
As stated earlier, we consider three model types introduced in \citet
{lauritzensym}, which can be represented by colored graphs. These are
discussed briefly here, and we refer to \citet{lauritzensym} for
further details.
\subsubsection*{RCON models: Restrictions on concentrations}
RCON models place\break \mbox{equality} constraints on the concentration matrix $K$.
They restrict off-diagonal elements of $K$ separately from those on the
diagonal, so that the restrictions can be represented by a graph
coloring $(\mathcal{V},\mathcal{E})$ of $G$, with $\mathcal{V}$
representing the diagonal and $\mathcal{E}$ the off-diagonal
constraints. The corresponding set of positive definite matrices is
denoted by $\mathcal{S}^{+}(\mathcal{V},\mathcal{E})$.

\subsubsection*{RCOR models: Restrictions on partial correlations}
RCOR models combine \mbox{equality} restrictions on the diagonal of $K$ with
equality constraints on partial correlations, given by
%
\begin{equation}\label{rho}
\rho_{\alpha\beta|V\setminus\{\alpha,\beta\}} = - \frac
{k_{\alpha\beta}}{\sqrt{k_{\alpha\alpha}k_{\beta\beta}}},\qquad \alpha
,\beta\in V, \alpha\not=\beta.
\end{equation}
%
Constraints of RCOR models may also be represented by a graph coloring
$(\mathcal{V},\mathcal{E})$, and we denote the corresponding set of
positive definite matrices by $\mathcal{R}^{+}(\mathcal{V},\mathcal{E})$.

\subsubsection*{RCOP models: Permutation symmetry}
RCOP models are determined by distribution invariance under a group of
permutations of the vertices which preserve the edges of the graph,
that is, a subgroup of $\operatorname{Aut}(G)$, the group of \emph
{automorphisms} of $G$. For $\sigma\in\operatorname{Aut}(G)$, let
$G(\sigma)$ be the permutation matrix representing $\sigma$, with
$G(\sigma)_{\alpha\beta} = 1$ if and only if $\sigma$ maps $\beta$
to $\alpha$, for $\alpha, \beta\in V$. Then a Gaussian $\mathcal
{N}_{|V|}(0,\Sigma)$ distribution is preserved by $\sigma$ if and
only if
%
\begin{equation}\label{cond}
G(\sigma)K G(\sigma)^{-1}= K.
\end{equation}
The RCOP model generated by a group $\Gamma\subseteq\operatorname
{Aut}(G)$ assumes
\[
K \in\mathcal{S}^{+}(G,\Gamma) = \mathcal{S}^{+}(G) \cap \mathcal
{S}^{+}(\Gamma),
\]
where $\mathcal{S}^{+}(\Gamma)$ denotes the set of positive definite
matrices satisfying (\ref{cond}) for all $\sigma\in\Gamma$ [\citet
{lauritzensym}].
\subsection{Relations between model types}
Under certain conditions on the coloring, RCON and RCOR models, which
are determined by the same colored graph, coincide in their model
restrictions. First we define edge regularity of a graph 
coloring.\vadjust{\goodbreak}

\begin{definition}\label{def_edge_reg}
Let $\mathcal{G}=(\mathcal{V},\mathcal{E})$ be a colored
graph. We say that $(\mathcal{V},\mathcal{E})$ is \emph{edge
regular} if any pair of edges in the same color class in $\mathcal{E}$
connects the same vertex color classes.
\end{definition}

The relevant result in \citet{lauritzensym} then becomes:

\begin{proposition}\label{prop_rcon_rcor}
The RCON and RCOR models, that are determined by the colored graph
$\mathcal{G}=(\mathcal{V},\mathcal{E})$, yield identical
restrictions, that is,
\[
\mathcal{S}^{+}(\mathcal{V},\mathcal{E}) = \mathcal{R}^{+}(\mathcal
{V},\mathcal{E})
\]
if and only if $(\mathcal{V},\mathcal{E})$ is edge regular.
\end{proposition}

RCOP models fall into the class of models which satisfy the condition
in Proposition~\ref{prop_rcon_rcor}, as we show below:

\begin{proposition}\label{rcop_edge}
If a colored graph $\mathcal{G}=(\mathcal{V},\mathcal{E})$
represents an RCOP model, then $(\mathcal{V},\mathcal{E})$ is edge regular.
\end{proposition}

\begin{pf} Let $\mathcal{G}=(\mathcal{V},\mathcal{E})$ represent
an RCOP model, generated by permutation group $\Gamma\subseteq
\operatorname{Aut}(G)$, say, and let $u \in\mathcal{E}$ and $e,f \in
u$. By definition, there exists $\sigma\in\Gamma$ mapping $e$ to $f$
while leaving $(\mathcal{V},\mathcal{E})$ invariant. This implies
that the colorings of the end vertices of $e$ and $f$ must be
identical.
\end{pf}

Thus, if for a graph $G=(V,E)$ we let $\mathcal{V}$ and $\mathcal{E}$
denote the vertex orbits and edge orbits of a group $\Gamma\subseteq
\operatorname{Aut} (G)$, then $\mathcal{S}^{+}(G,\Gamma) = \mathcal
{S}^{+}(\mathcal{V},\mathcal{E}) = \mathcal{R}^{+}(\mathcal
{V},\mathcal{E})$ [\citet{lauritzensym}].

For example, since the coloring of the graph in Figure~\ref{fig:frets}
is generated by the group $\Gamma=\{I, \sigma\}$ with $\sigma$
simultaneously permuting $B_{1}$ with $B_{2}$ and $L_{1}$ with $L_{2}$,
the corresponding two sets $\mathcal{S}^{+}(\mathcal{V},\mathcal
{E})$ and $\mathcal{R}^{+}(\mathcal{V},\mathcal{E})$ coincide.

\section{Equality of maximum likelihood and least-squares estimator}
By Theorem~\ref{thm:kruskal}, the maximum likelihood estimator $\hat
{\mu}$ and least-squares estimator $\mu^{*}$ agree for $\mu$ in a
linear subspace $\Omega\subseteq\mathbb{R}^{V}$ if and only if
$\Omega$ is stable under all $K$ in the model. Below we show that for
RCON and RCOR models, and thus also for RCOP models, invariance of
$\Omega$ under $K$ is equivalent to invariance under $\{T^{u}\}_{u \in
\mathcal{V} \cup\mathcal{E}}$.

\begin{proposition}\label{invar_rcon}
Let $\mathcal{G} = (\mathcal{V},\mathcal{E})$ be a colored graph
representing the RCON model with $K \in\mathcal{S}^{+}(\mathcal
{V},\mathcal{E})$. Then for $\Omega\subseteq\mathbb{R}^{V}$,
\[
K\Omega\subseteq\Omega \qquad\forall K \in\mathcal{S}^{+}(\mathcal
{V},\mathcal{E}) \quad\Longleftrightarrow\quad T^{u} \Omega\subseteq\Omega
\qquad\forall u \in\mathcal{V} \cup\mathcal{E}.
\]
\end{proposition}
\begin{pf} By definition of RCON models, all $K \in\mathcal
{S}^{+}(\mathcal{V},\mathcal{E})$ can be written as
%
\begin{equation}\label{k_restrictions}
K = \sum_{u \in\mathcal{V} \cup\mathcal{E}} \theta_{u}T^{u},\qquad
\theta_{u} \in\mathbb{R} \mbox{ for } u \in\mathcal{V} \cup
\mathcal{E}\vadjust{\goodbreak}
\end{equation}
with $\{\theta_{u}\}_{u \in\mathcal{V} \cup\mathcal{E}}$ such that
the expression in (\ref{k_restrictions}) is positive definite.
Suppose first that we have
\[
T^{u} \Omega\subseteq\Omega \qquad\forall u \in\mathcal{V} \cup\mathcal{E}.
\]
Since $\Omega$ is a linear space this implies invariance under all $K
\in\mathcal{S}^{+}(\mathcal{V},\mathcal{E})$.

Next suppose $K\Omega\subseteq\Omega$ for all $K \in\mathcal
{S}^{+}(\mathcal{V},\mathcal{E})$. $\mathcal{S}^{+}(\mathcal
{V},\mathcal{E})$ is an open convex cone, so that for all $K \in
\mathcal{S}^{+}(\mathcal{V},\mathcal{E})$ and $u \in\mathcal{V}
\cup\mathcal{E}$ there exists $\lambda_{u} \in\mathbb{R} \setminus
\{0\}$ such that
\[
K_{u} = K + \lambda_{u} T^{u} \in\mathcal{S}^{+}(\mathcal
{V},\mathcal{E}).
\]
By assumption, $K\Omega\subseteq\Omega$ and $K_{u}\Omega\subseteq
\Omega$, which gives
\[
(K_{u}-K)\Omega= \lambda_{u} T^{u}\Omega\subseteq\Omega
\]
and thus the desired result.
\end{pf}

Although the cone $ \mathcal{R}^{+}(\mathcal{V},\mathcal{E})$ is not
in general convex, the same holds for RCOR models, as shown below.
\begin{proposition}\label{invar_rcor}
Let $\mathcal{G}=(\mathcal{V},\mathcal{E})$ be a colored graph. Then
for $\Omega\subseteq\mathbb{R}^{V}$,
\[
K\Omega\subseteq\Omega \qquad\forall K \in\mathcal{R}^{+}(\mathcal
{V},\mathcal{E}) \quad\Longleftrightarrow\quad T^{u} \Omega\subseteq\Omega
\qquad\forall u \in\mathcal{V} \cup\mathcal{E}.
\]
\end{proposition}

\begin{pf} Let $A=(a_{\alpha\beta})_{\alpha, \beta\in V}$
denote the diagonal matrix with entries equal to the inverse partial
standard deviations, that is,
\[
a_{\alpha\alpha} = \sqrt{k_{\alpha\alpha}}\qquad \mbox{for } \alpha
\in V,
\]
and let $C=(c_{\alpha\beta})_{\alpha, \beta\in V}$ have all
diagonal entries equal to 1 and all off-diagonal entries be given by
the negative partial correlations $- \rho_{\alpha\beta| V\setminus\{
\alpha, \beta\}}$. Then, by (\ref{rho}), all $K \in\mathcal
{R}^{+}(\mathcal{V},\mathcal{E})$ can be uniquely expressed as
%
\begin{equation}\label{k_repar}
K = ACA = \sum_{v,w \in\mathcal{V}} a_{v}a_{w}T^{v}T^{w} + \sum
_{v,w \in\mathcal{V}, u \in\mathcal{E}} c_{u}a_{v}a_{w}T^{v}T^{u}T^{w},
\end{equation}
where for $v \in\mathcal{V}$ and $u \in\mathcal{E}$, we let $a_{v}$
and $c_{u}$ denote $a_{\alpha\alpha}$, $\alpha\in v$ and $c_{\alpha
\beta}$, $\alpha\beta\in u$, respectively. As for $v,w \in\mathcal
{V}$, $T^{v}T^{w}$ is zero unless $v=w$, when it equals $T^{v}$, and
equation (\ref{k_repar}) simplifies to
%
\begin{equation}\label{eq:k_in_rcor_expr}
K =\sum_{v \in\mathcal{V}} a_{v}^{2}T^{v} + \sum_{v, w \in\mathcal
{V}, u \in\mathcal{E}} c_{u}a_{v}a_{w}T^{v}T^{u}T^{w}.
\end{equation}
Suppose first that we have
\[
T^{u} \Omega\subseteq\Omega\qquad \forall u \in\mathcal{V} \cup\mathcal{E}.
\]
As before, since $\Omega$ is a linear space, (\ref
{eq:k_in_rcor_expr}) implies invariance under all $K \in\mathcal
{R}^{+}(\mathcal{V},\mathcal{E})$. So suppose $K\Omega\subseteq
\Omega$ for all\vadjust{\goodbreak} $K \in\mathcal{R}^{+}(\mathcal{V},\mathcal{E})$.
Then, in particular, $\Omega$ is invariant under all $K$ of the form
\[
K=ACA\qquad \mbox{with $A = \sigma^{2}I$ and $C$ as specified above},
\]
represented by the graph coloring $(\{V\},\mathcal{E})$, which has the
same edge coloring as $(\mathcal{V},\mathcal{E})$, but with all
vertices of the same color. This graph coloring is clearly edge regular
(as all edges connect the only vertex color class back to itself),
which gives that the represented model is also of type RCON. We may
therefore apply Proposition~\ref{invar_rcon}, giving
\[
T^{u} \Omega\subseteq\Omega \qquad\forall u \in\mathcal{E}.
\]

For the vertex coloring, consider the submodel
\[
K=ACA \qquad\mbox{with $A$ as specified above and $C=I$}.
\]
This submodel is represented by the independence graph with no edges
and vertex coloring $\mathcal{V}$, which is also edge regular, so that
by Proposition~\ref{invar_rcon},
\[
T^{v} \Omega\subseteq\Omega\qquad \forall v \in\mathcal{V},
\]
completing the proof.
\end{pf}

\subsection{Equality restrictions in the means of RCON and RCOR models}
Let~$\mathcal{M}$ be a partition of $V$ and consider $\Omega(\mathcal
{M})$ as in (\ref{omega_m}).
In the following we derive a necessary and sufficient condition on
$\mathcal{M}$ and $(\mathcal{V},\mathcal{E})$ for $\Omega(\mathcal
{M})$ to be invariant under all $K \in\mathcal{S}^{+}(\mathcal
{V},\mathcal{E})$ or all $K \in\mathcal{R}^{+}(\mathcal{V},\mathcal
{E})$. By Propositions~\ref{invar_rcon} and~\ref{invar_rcor} we may
consider vertex and edge color classes of the colored graph
representing the model separately. We begin with the vertex coloring.

\begin{proposition}\label{prop_finer_par}
Let $\mathcal{M}$ and $\mathcal{V}$ be partitions of $V$. Then
$\Omega(\mathcal{M})$ is invariant under $\{T^{v}\}_{v \in\mathcal
{V}}$ if and only if $\mathcal{M} \leq\mathcal{V}$.
\end{proposition}

\begin{pf} The action of $T^{v}$ for $v \in\mathcal{V}$ on $\mu
\in\mathbb{R}^{V}$ is given as
\[
T^{v}\mu= \cases{
\mu_{\alpha} &\quad if $\alpha\in v$,
\cr
0 &\quad otherwise.}
\]
Let $\mu\in\Omega(\mathcal{M})$, and suppose that $T^{v}\Omega
(\mathcal{M}) \subseteq\Omega(\mathcal{M})$ for all $v \in\mathcal
{V}$. In order for $T^{v}\mu\in\Omega(\mathcal{M})$ for all $v \in
\mathcal{V}$, we must have $\alpha\equiv\beta (\mathcal{V})$
whenever $\alpha\equiv\beta (\mathcal{M})$, or equivalently
$\mathcal{M} \leq\mathcal{V}$. Conversely, suppose $\mathcal{M}
\leq\mathcal{V}$. Then $\alpha\equiv\beta (\mathcal{V})$ whenever
$\alpha\equiv\beta (\mathcal{M})$, which gives $T^{v}\mu\in\Omega
(\mathcal{M})$ for all $v \in\mathcal{V}$.
\end{pf}

Note that the above result implies that the likelihood cannot be
maximized in $\mu$ independently of the value of $K$ in the
Behrens--Fisher setting, which is the RCON and RCOR model on two
variables specified by $\mathcal{V}=\{\{1\},\{2\}\}$, $\mathcal
{E}=\varnothing$ together with the restriction $\mu_{1} = \mu_{2}$ on
the means, as the mean partitioning is then coarser than the vertex coloring.

Also, for the model of Frets's heads in Figure~\ref{fig:frets}, the
MLE of the mean is not simple if we, for example, wish to estimate the
mean under the hypothesis that the heads tend to be square-shaped, that
is, if mean lengths are equal to mean breadths, as this partition would
not be finer than the vertex coloring.

For the edge coloring, we require the concept of an \emph{equitable
partition}, first defined by \citet{sachs}.

\begin{definition}[(Sachs)]\label{equitable_def}
Let $G=(V,E)$ be an undirected graph. Then a~partition, or
equivalent coloring, $\mathcal{V}$ of $V$ is called \emph{equitable
with respect to~$G$} if for all $v \in\mathcal{V}$, $\alpha, \beta
\in v$,
\[
|\operatorname{ne}_{E}(\alpha) \cap w| = |\operatorname
{ne}_{E}(\beta) \cap w| \qquad\forall w \in\mathcal{V}.
\]
\end{definition}

\citet{eigen} proved the following. 

\begin{proposition}[(Chan and Godsil)]\label{prop_chan}
Let $G=(V,E)$ be an undirected graph with adjacency matrix $T$, and let
$\mathcal{M}$ be a partition of $V$. Then $\Omega(\mathcal{M})$ is
invariant under $T$ if and only if $\mathcal{M}$ is equitable with
respect to $G$.
\end{proposition}

The notion of an equitable partition for vertex colored graphs can be
naturally extended to graphs with colored vertices and edges. We term
the corresponding graph colorings \emph{vertex regular}, defined below.

\begin{definition}\label{palette}
Let $\mathcal{G}=(\mathcal{V},\mathcal{E})$ be a colored
graph, and let the subgraph induced by the edge color class $u \in
\mathcal{E}$ be denoted by $G^{u}=(V,u)$. We say that $(\mathcal
{V},\mathcal{E})$ is \emph{vertex regular} if $\mathcal{V}$ is
equitable with respect to $G^{u}$ for all $u \in\mathcal{E}$.
\end{definition}

Combining Definition~\ref{palette} with Proposition~\ref{prop_chan}
yields vertex regularity to be a~necessary and sufficient condition for
$\Omega(\mathcal{M})$ to be invariant under $\{T^{u}\}_{u \in
\mathcal{E}}$:

\begin{proposition}\label{prop_edge_invariance}
For $\mathcal{G}=(\mathcal{V},\mathcal{E})$ a colored graph and
$\mathcal{M}$ a partition of~$V$, $\Omega(\mathcal{M})$ is invariant
under $\{T^{u}\}_{u \in\mathcal{E}}$ if and only if $(\mathcal
{M},\mathcal{E})$ is vertex regular.
\end{proposition}

\begin{pf} For $u \in\mathcal{E}$, $T^{u}$ is the adjacency
matrix of $G^{u}$. By Proposition~\ref{prop_chan}, $\Omega(\mathcal
{M})$ is stable under $\{T^{u}\}_{u \in\mathcal{E}}$ if and only if
$\mathcal{M}$ is equitable with respect to $G^{u}$ for all $u \in
\mathcal{E}$, or equivalently if and only if $(\mathcal{M},\mathcal
{E})$ is vertex regular.
\end{pf}

For the Frets's heads model in Figure~\ref{fig:frets}, this implies
that restricting the mean breadths and mean lengths on their own does
not ensure $\hat{\mu}=\mu^{*}$, as the corresponding partitions
$\mathcal{M}=\{\{B_{1},B_{2}\},\{L_{1}\},\{L_{2}\}\}$ and $\mathcal
{M}=\{\{B_{1}\},\{B_{2}\},\{L_{1},L_{2}\}\}$ do not give rise to vertex
regular colorings $(\mathcal{M},\mathcal{E})$.

Combining Proposition~\ref{prop_finer_par} and Proposition~\ref
{prop_edge_invariance} establishes our main result:

\begin{them}\label{means}
Let $\mathcal{G}=(\mathcal{V},\mathcal{E})$ be a colored graph, let
$\mathcal{M}$ be a partition of $V$ and consider a sample\vadjust{\goodbreak} from a
multivariate Gaussian $\mathcal{N}_{|V|}(\mu,\Sigma)$ distribution
with $\mu\in\Omega(\mathcal{M})$ and $K = \Sigma^{-1} \in\mathcal
{S}^+(\mathcal{V},\mathcal{E})$, both unknown. It then holds that
%
\begin{equation}\label{thm_cond}
\hat{\mu} = \mu^{*} \quad\Longleftrightarrow\quad \mathcal{M} \leq\mathcal
{V}\quad \mbox{and}\quad (\mathcal{M},\mathcal{E}) \mbox{ is vertex regular.}
\end{equation}
The same conclusion holds if $\mathcal{S}^+(\mathcal{V},\mathcal
{E})$ is replaced by $\mathcal{R}^+(\mathcal{V},\mathcal{E})$.
\end{them}

In fact, for any colored graph $\mathcal{G}=(\mathcal{V},\mathcal
{E})$ there is always a coarsest partition $\mathcal{M}$ satisfying
the conditions in (\ref{thm_cond}) [\citet{gehrmann}]. The finest
variant of such a coarsest equitable refinement $\mathcal{M}$ of
$\mathcal{V}$ is given by the singleton partition. Clearly it is finer
than any vertex coloring, and further naturally gives a vertex regular
coloring $(\mathcal{M},\mathcal{E})$. Note that in this case $\Omega
(\mathcal{M})$ corresponds to the unrestricted case considered in
Section~\ref{mle_section}, conforming with the fact that then $\hat
{\mu} = \mu^{*}$.

\subsection{Equality restrictions in the means in RCOP models}
The coarsest possible $\mathcal{M}$ for RCON and RCOR models, by (\ref
{thm_cond}), is $\mathcal{V}$, for which $\hat{\mu} = \mu^{*}$ if
and only if $(\mathcal{V},\mathcal{E})$ is vertex regular. This
always holds for RCOP models.

\begin{proposition}\label{rcop_palette}
If a colored graph $\mathcal{G}=(\mathcal{V},\mathcal{E})$
represents an RCOP model, then $(\mathcal{V},\mathcal{E})$ is vertex regular.
\end{proposition}

\begin{pf} Let $\mathcal{G}=(\mathcal{V},\mathcal{E})$ represent
an RCOP model, generated by a subgroup $\Gamma\subseteq\operatorname
{Aut}(G)$, say. By Proposition~\ref{rcop_edge}, $(\mathcal
{V},\mathcal{E})$ is edge regular. Thus whenever two edges $e,f \in E$
are of the same color, they connect the same vertex color classes.
Let $\alpha, \beta\in V$ be two equally colored vertices in $\mathcal
{V}$. Then, by definition of RCOP models, there exists a permutation
$\sigma\in\Gamma$ which maps $\alpha$ to $\beta$ leaving
$(\mathcal{V},\mathcal{E})$ invariant. This implies that the degree
in each edge color class of $\alpha$ and $\beta$ must be identical.
The previous two statements imply
\[
|\operatorname{ne}_{u}(\alpha) \cap v| = |\operatorname
{ne}_{u}(\beta) \cap v|
\]
for all $v \in\mathcal{V}$ and all pairs $\alpha, \beta\in V$ with
$\alpha\equiv\beta (\mathcal{V})$, which is precisely the criterion
of vertex regularity for the graph coloring $(\mathcal{V},\mathcal
{E})$.
\end{pf}

We mention in passing that Proposition~\ref{rcop_edge} and
Proposition~\ref{rcop_palette} combined establish that colorings of
graphs which represent RCOP models are \emph{regular} in the
terminology of \citet{siemonsaut}.
We conclude from Theorem~\ref{means} and Proposition~\ref{rcop_palette}:

\begin{corollary}\label{cor}
Let $G=(V,E)$ be an undirected graph, and let $(\mathcal{V},\mathcal
{E})$ represent the constraints of an RCOP model generated by group
$\Gamma\subseteq\operatorname{Aut}(G)$. Then for a sample from a
multivariate Gaussian $\mathcal{N}_{|V|}(\mu,\Sigma)$ distribution
with $\mu\in\Omega(\mathcal{V})$ and $K = \Sigma^{-1} \in\mathcal
{S}^{+}(G,\Gamma)$, both unknown, we always have $\hat{\mu} = \mu^{*}$.
\end{corollary}

\section{Examples} \label{examples}

We first consider the example in Figure~\ref{fig:frets} on head
dimensions for first and second sons. Representing an RCOP model, by
Proposition~\ref{rcop_palette}, the graph in
Figure~\ref{fig:frets}
has a vertex regular coloring. It follows from Corollary~\ref{cor}\vadjust{\goodbreak}
that the maximum likelihood estimate of the mean under the hypothesis
that the mean length and mean breadth are equal for the two sons is
simply the total average of the head lengths and the head breadths,
respectively. Similarly, it follows from Theorem~\ref{means} that the
only hypotheses about the mean that have a simple solution are this one
and the one where the means are completely unrestricted. 

The empirical means of the dimensions $(B_1,B_2,L_1,L_2)$ are equal to
$(151.12,\allowbreak 149.24, 185.72, 183.84)$, so that the MLE of the means under
the hypothesis that the mean lengths and breadths are independent of
the parity of the son then become $(150.18,150.18,184.78, 184.78)$. The
likelihood ratio test is obtained by comparing the maximized profile
likelihoods (\ref{eq:proflik}) calculated with appropriate residual
covariance matrices $W$ under the two hypotheses. Using the \texttt
{R}-package \texttt{gRc} [\citet{grc}] this yields $-2\log\mathrm
{LR}=3.27$ on 2 degrees of freedom, so there is no evidence for the
sizes depending on the parity of the son.

\begin{figure}

\includegraphics{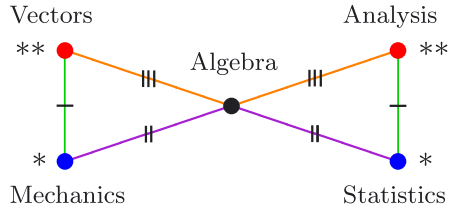}

\caption{Mathematics marks example.}
\label{fig_math}
\end{figure}

Our second example is concerned with the examination marks of 88
students in five mathematical subjects [\citet{mardia}]. The RCOP model
represented by $\mathcal{G}=(\mathcal{V},\mathcal{E})$ in
Figure~\ref{fig_math} was demonstrated to be an excellent fit in \citet
{lauritzensym}.

%

The model is given by invariance of $K$ under simultaneously permuting
Mechanics with Statistics, and Vectors with Analysis. In their fit,
\citet{lauritzensym} implicitly assumed an unconstrained mean. The
MLE~$\hat{\mu}$ is then given by the sample averages $(\hat{\mu
}_{\mathrm{al}}, \hat{\mu}_{\mathrm{an}}, \hat{\mu}_{\mathrm
{me}}, \hat{\mu}_{\mathrm{st}}, \hat{\mu}_{\mathrm{ve}}) =
(50.60,46.68,\allowbreak 38.96,42.31, 50.59)$ in the obvious notation, which
corresponds to $\mathcal{M}$ being the singleton partition.
However, it could be natural to assume $\mu$ subject to the same
invariance as $K$, meaning $\mathcal{M} = \mathcal{V}$, or in this
case $\mu_{\mathrm{an}}=\mu_{\mathrm{ve}}$ and $\mu_{\mathrm
{me}}=\mu_{\mathrm{st}}$. Then, by Corollary~\ref{cor}, $(\hat{\mu
}_{\mathrm{al}}, \hat{\mu}_{\mathrm{an}}, \hat{\mu}_{\mathrm
{me}}, \hat{\mu}_{\mathrm{st}}, \hat{\mu}_{\mathrm{ve}}) =
(50.60,48.64,40.63,\allowbreak 40.63,48.64)$. The likelihood ratio statistic for
this mean structure relative to the model in Figure~\ref{fig_math}
with unconstrained mean takes the value of 11.9 on 2 degrees of
freedom, and the hypothesis about symmetry in the means is therefore
clearly rejected with $p<0.003$.

Note that the sample averages for Vectors and Algebra are almost
identical. However, combining the constraints on $K$ represented by the
graph in Figure~\ref{fig_math} with any hypothesis on the means
implying $\mu_{\mathrm{ve}} = \mu_{\mathrm{al}}$ will require joint
maximization in $\mu$ and $K$ of the likelihood function in (\ref
{eq:likfunction}) to obtain $\hat{\mu}$.

\section{Discussion}

The main result of this article is a necessary and sufficient condition
on the pattern of equality constraints on the mean vector $\mu$ in a~graphical Gaussian symmetry model with colored graph $\mathcal
{G}=(\mathcal{V},\mathcal{E})$ which ensures the identity of the
least squares estimate $\mu^{*}$ and the maximum likelihood estimate
$\hat{\mu}$, given in Theorem~\ref{means}.

The derived necessary and sufficient condition is formulated in terms
of vertex and edge colored graphs and is easily testable, so that any
set of equality constraints on $\mu$ together with constraints on
either $K$ or the partial correlations can be verified for estimability
of $\mu$ by $\mu^{*}$.

The result is, for example, useful if one set of constraints, either on
the mean or the independence structure, is assumed to be given and the
other may be varied. A setting which falls into this category is the
design of experiments seeking to estimate mean treatment effects under
the assumption of correlations with inherent symmetries in the error
structure at experimental sites. A systematic arrangement of the sites
may enforce a symmetry pattern in the concentrations or partial
correlations, and thus restrict the concentration matrix as in one of
models considered here. Allocation of treatments then effectively
restricts the mean response at sites with the same treatment to be
identical, and the condition derived can be used to find treatment
allocations which ensure estimability of mean treatment effects without
knowledge of the value of $\Sigma$.

An interesting question directly emerging from our work is concerned
with the exact distributions of likelihood ratio test statistics for
hypotheses about the mean in RCON, RCOR and RCOP models. In the
examples discussed in this paper we have relied on asymptotic theory to
judge the significance of test statistics, but it is likely that their
distributions can be derived explicitly, for example when the mean
hypothesis is given by the natural symmetry of an RCOP model. \citet
{hylleberg} developed explicit likelihood ratio tests for decomposable
mean-zero RCOP models generated by compound symmetry [\citet{votaw}],
and it would be interesting to extend these results to models with
nonzero means and more general symmetry constraints.

We note that if our condition is not satisfied, we would expect
phenomena similar to those in the Behrens--Fisher problem implying the
general nonuniqueness of the MLE [\citet{drtonbehrens}] and the
nonexistence of an $\alpha$-similar test of the mean hypothesis
[\citet{scheffe}].


\printaddresses

\end{document}